\documentclass{commat}

\newcommand{\DD}{\mathbf{D}}
\newcommand{\A}{\mathcal{A}}
\newcommand{\Q}{\mathcal{Q}}
\newcommand{\QQ}{\mathbf{Q}}
\title{%
    A complete invariant for doodles on a 2-sphere
    }

\author{%
    Jacob Mostovoy
    }

\authorinfo[%
    Jacob Mostovoy]{% Please, use "F. Name" format
    CINVESTAV, Mexico City, Mexico}{%
    jacob@math.cinvestav.mx
    }

\abstract{%
    We define a complete invariant for doodles on a 2-sphere which takes values in a series of chord diagrams of a certain type. The coefficients in the diagrams with $n$ chords are finite type invariants of doodles of order at most $2n$.
    }



\VOLUME{31}
\YEAR{2023}
\NUMBER{3}
\firstpage{137}
\DOI{https://doi.org/10.46298/cm.12893}

\begin{document}

%\begin{center}
 {\textit{{To the memory of Sergei Duzhin}}}

%\includegraphics[width=400 pt]{Duzhin.jpg}

%\end{center}
\section{Introduction}
    The theory of finite-type knot invariants may be thought of as 
a nilpotency theory for knots. Nilpotent and residually nilpotent groups can be studied by embedding them into algebras of power series; most notably, the free group on $n$ letters $x_1, \ldots, x_n$ maps into the algebra of the noncommutative power series in $X_1, \ldots, X_n$  by
$$x_i\mapsto 1+ X_i.$$
One may describe this embedding by saying that a word that involves only non-negative exponents is mapped to the formal sum of its subwords \cite{MKS}.

In a similar fashion, the set of isotopy classes of knots maps to the power series in chord diagrams by means of the Kontsevich integral. In this case, however, the injectivity of this map is a major open problem. Moreover, the combinatorial interpretation of the Kontsevich integral of a knot is rather sophisticated and involves Drinfeld associators \cite{CDBook}. 

The Kontsevich integral contains complete information about the finite type invariants of a knot. One may be tempted to construct universal finite-type knot invariants by using the same elementary method that works for the free groups. Namely, one can encode knots by diagrams of some kind and then assign to a knot a formal sum of the subdiagrams of its diagram. In particular, if knots are represented by so-called clasp diagrams  \cite{MP}, the formal sum of the subdiagrams taken by modulo diagrams with $n+1$ clasps is, in fact, a universal invariant of order $n$. This construction, however, is no replacement for the Kontsevich integral since the values of this invariant lie in a certain algebra which is rather more involved than the (already uncomputable) algebra of chord diagrams. Essentially, it is a restatement of the original approach of Vassiliev in different terms.

The combinatorics of knots is much more complex than that of words in a free group. However, for other knotted objects, the approach via sums of subdiagrams may be very efficient. Pure braids, which are closely related to knots, form groups that are almost direct products of free groups. Each pure braid has a unique normal (``combed'') form; the formal sum of the subwords of such a normal form is a universal finite type invariant whose values lie in the same algebra as the Kohno integral for braids, of which the Kontsevich integral is a generalization. This universal invariant is injective and this is one way to show that the finite type invariants distinguish non-equivalent pure braids \cite{MW, Pap}.

In the present note, we consider the set of equivalence classes of curves in $S^2$ without triple self-intersections, called ``doodles'' by Khovanov\footnote{The term ``doodles'' was first used \cite{FT} for a different class of curves. Khovanov's terminology has now superseded that of Fenn and Taylor.} in \cite{Kh}. We show that the sums of subdiagrams give a  finite type invariant of curves of doodles in $S^2$, which is complete in the sense that it distinguishes non-equivalent curves.

The theory of curves in $S^2$ without triple self-intersections may be thought of as a toy version of knot theory. Such curves are encoded by Gauss-type diagrams; the equivalence is generated by certain moves similar to the Reidemeister moves.  These curves may be presented as closures of braid-like groups \cite{Kh, Got} and possess invariants similar to the knot group and to polynomial invariants \cite{Juy}. The radical difference between doodles and knots in $S^3$ lies in the fact that the problem of distinguishing doodles is, essentially, trivial: there exists a simplification algorithm that produces a unique representative for each equivalence class of doodles, and this representative is minimal in the number of crossings. In view of this circumstance, the theory of doodles is mostly important as a testing ground for knot theoretic constructions.

The finite-type invariants of doodles have been studied by Vassiliev and Merkov; not only in $S^2$ but also in the more difficult case of doodles in ${\mathbb R}^2$, see \cite{Vas1, Merk}. In particular, among other things, Merkov showed that finite type invariants of doodles in $S^2$ distinguish non-equivalent doodles. His method is based on the use of the Polyak-Viro type invariants that count the subdiagrams of the arrow diagram of a doodle. What we do is very close in spirit and is based on the work of Goussarov, Polyak and Viro \cite{GPV}, who defined the finite type invariants for virtual knots in terms of a certain universal invariant\footnote{Their definition is different from the definition of the finite type invariants for virtual knots given by Kauffman in \cite{Ka}.}. In the case of knots, the Polyak-Viro type knot invariants span the space of all finite type invariants (this result is known as the Goussarov Theorem) while the Goussarov-Polyak-Viro universal invariant for ``compact'' virtual knots is strictly weaker. As we shall see, in the case of doodles in $S^2$ both approaches produce complete sets of invariants and, therefore, are equivalent in strength.

\section{Doodles and their arrow diagrams}

A \emph{doodle} is an immersion $S^1\to S^2$ with a finite number of transversal self-intersections (or \emph{crossings}). Doodles are considered modulo isotopies of $S^2$ and the following two moves:   
$$\includegraphics[width=3in]{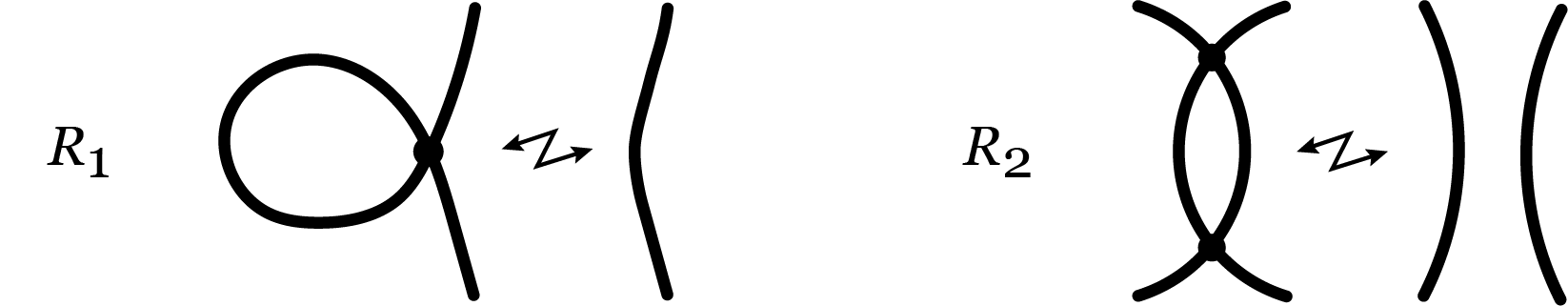}$$  
These moves can also be considered on curves in $\mathbb{R}^2$ or in other surfaces as it is done in \cite{BFKK}. 
We will refer to two doodles that can be connected by these transformations as being equivalent. A \emph{trivial} doodle is a doodle without crossing points. All trivial doodles in $S^2$ are equivalent.

A doodle in $S^2$ is determined up to isotopy by its \emph{arrow diagram}\footnote{It would be more natural to use the term ``Gauss diagram''. However, Arnold in \cite{Ar} uses the term ``Gauss diagram'' for the underlying chord diagram of the arrow diagram, so we follow the terminology of Merkov. The term ``arrow diagram'' has also been used in a somewhat different context in \cite{Pol}.}, which is a chord diagram with directed chords. The skeleton (that is, the outer circle) of the chord diagram corresponds to the parameter along the curve, the chords connect the inverse images of the crossings. At each crossing the branches of the doodle can be ordered in such a way that the tangent vectors to the first and the second branches form a positive basis of the tangent space to $S^2$. Each chord of the arrow diagram of the doodle is directed from the first branch to the second branch.

\begin{figure}[ht]
$$\includegraphics[height=1.2in]{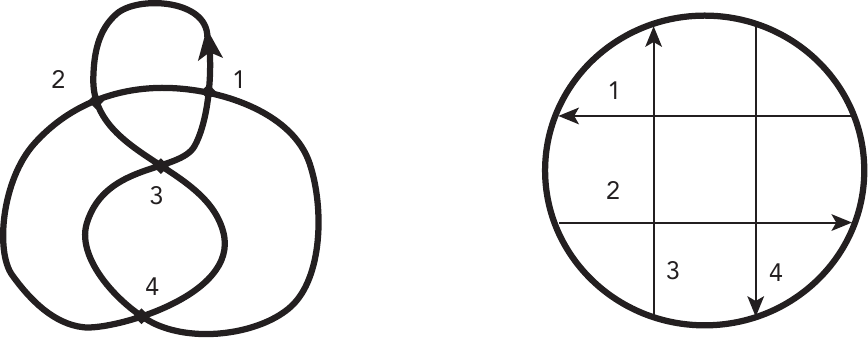}$$
\caption{A doodle and its arrow diagram}\label{arrow}
\end{figure}
 
The moves $R_1, R_2$ on doodles translate into the moves $R'_1, R'_2$ on their arrow diagrams: 
$$\includegraphics[height=0.6in]{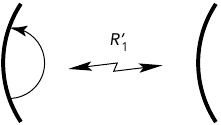}$$
and
$$\includegraphics[height=0.6in]{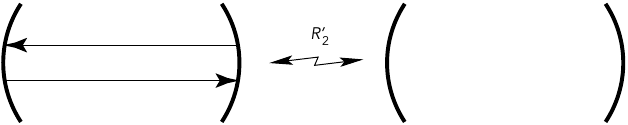}$$
Here, and in what follows, the convention is that we show only the relevant part of the arrow diagram; the omitted chords 
are the same in all the diagrams that appear in the same move or equation. Unless explicitly indicated, 
the orientations on the segments of the skeleton may be arbitrary but the same on all the diagrams. 

Not each arrow diagram comes from a doodle; if it does, it is called \emph{realizable}. 
The above moves, however, make sense for all arrow diagrams. We consider two arrow diagrams that are connected by these moves to be equivalent.

A doodle may be simplified by means of the moves $R_1, R_2$ that decrease the number of crossings. If a doodle cannot be simplified, it is called \emph{minimal}. Each doodle is equivalent to a unique (up to isotopy) minimal doodle \cite{Kh}. The same holds for arrow diagrams and the moves $R'_1, R'_2$; the proof of this fact is, essentially, the same as in the case of doodles.

Note that by applying a move $R'_1$ or $R'_2$ in a way that decreases the number of crossings to a realizable arrow diagram we again obtain a realizable diagram. It follows that the minimal diagram equivalent to a realizable arrow diagram is also realizable. We arrive at the following
\begin{proposition}\label{subset} If two realizable arrow diagrams are equivalent, the corresponding doodles are also equivalent.
\end{proposition}
In other words, the set of equivalence classes of doodles is a subset of the set of equivalence classes of arrow diagrams. In what follows we will identify doodles and realizable arrow diagrams.

A function on the set of isotopy classes of doodles that are preserved by the moves  $R_1, R_2$  is called an \emph{invariant} of doodles in $S^2$. Proposition~\ref{subset} implies that any function on the arrow diagrams that is preserved by the moves $R'_1, R'_2$  gives rise to an invariant of doodles and that each doodle invariant arises in this way.

Write $\DD$ for the set of all arrow diagrams and $\QQ\DD$ for the vector space over $\QQ$ generated by this set. The moves $R'_1, R'_2$ define a set of relations in $\QQ\DD$ which equate the left-hand side of a move to the right-hand side:
$$
\includegraphics[height=0.6in]{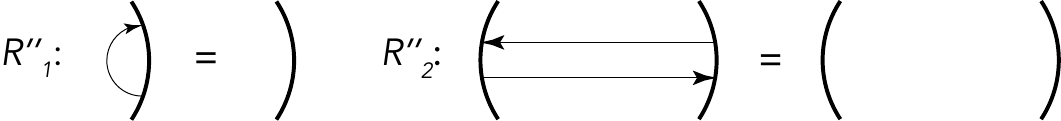}
$$
Define the linear map 
$$I_\DD: \QQ\DD\to \QQ\DD$$
by sending a diagram $D$ into the sum of all the diagrams obtained from $D$ by deleting a subset of chords of $D$.
The relations $R''_1$ and $R''_2$ in $\QQ\DD$ are sent by $I_\DD$ to the following relations:\begin{equation}\label{l1}
\includegraphics[height=0.6in]{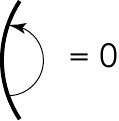}
\end{equation}
\smallskip
and
\smallskip
\begin{equation}\label{l2}
\includegraphics[height=0.5in]{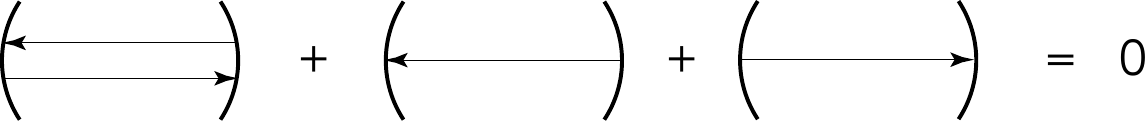}
\end{equation}
Denote by $\A$ the quotient of $\QQ\DD$ by the relations (\ref{l1}) and (\ref{l2}). The map $I_\DD$ composed with the quotient map gives rise to the map
$$I: \QQ\DD\to \A,$$
which, by construction, is an invariant of doodles in $S^2$.

For $n\geq 0$ define $\A_n$ to be the quotient of $\A$ by the subspace spanned by the diagrams with more than $n$ chords. There is a natural projection $\A_{n+1}\to \A_n$ which annihilates the diagrams with $n+1$ chords.
Write $I_n$ for the composition 
$$\QQ\DD\xrightarrow{I} \A \to \A_n.$$
The main result of the present note is the following
\begin{theorem}\label{main1}
The map $$I_n:\QQ\DD\to \A_n$$
gives rise to a finite type invariant of doodles of order at most $2n$. For a non-trivial doodle $D$ with $n$ or fewer crossings, 
$$I_n(D)\neq 0.$$ 
If $D_1, D_2$ are non-equivalent doodles with $n$ or fewer crossings,  
$$I_{n+1}(D_1)\neq I_{n+1}(D_2).$$
\end{theorem} 
We will refer to the doodle invariant defined by $I_n$ (which we will also denote by $I_n$) as the \emph{$n$-th diagram invariant}.  

For each $n$, the invariant $I_n$ is the composition of $I_{n+1}$ with the projection of  $\A_{n+1}$ onto $\A_n$. Let $\overline{\A}$ be the inverse limit of the spaces $\A_i$ and  $$\overline{I}: \QQ\DD\to \overline{\A}$$ the inverse limit of the $I_i$. Theorem~\ref{main1} implies
\begin{corollary}
The map $\overline{I}$ is a complete invariant of doodles in $S^2$.
\end{corollary} 

Let $V_n$ be the kernel of the projection $\A_{n+1}\to\A_n$. The space $\overline{\A}$ is isomorphic to the product of all the $V_n$. Another consequence of Theorem~\ref{main1} is 
\begin{corollary}
The map $\overline{I}$ composed with the projection of $\overline{\A}$ to $V_n$ is a finite type invariant of order at most $2n$.
\end{corollary} 
If a basis is chosen in each $V_n$, the map $\overline{I}$ becomes a series of arrow diagrams whose coefficients are finite type invariants, and the coefficient at any basis diagram in $V_n$ is a finite type invariant of order at most $2n$.

In fact, we will see that Theorem~\ref{main1}  holds even if $\A_n$ is replaced by a certain quotient of $\A_n$ which admits an explicit description. We consider it in the next section. In Section~\ref{Vassiliev} we will show that $I_n$ is an invariant of order $2n$.

\section{Quiver diagrams and the diagram invariants} 

Let us say that two chords in a chord diagram are \emph{adjacent} if their endpoints are adjacent:
$$\includegraphics[width=1.4in]{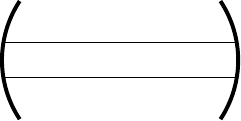}$$

If a chord diagram $C$ has a distinguished pair of adjacent chords, $C$ can be \emph{reduced} by deleting one of the chords in the pair. The reduction of a diagram $C$ defines a map from the set of chords of $C$ to the set of chords of the reduced diagram $C'$: it maps both adjacent chords to one chord that replaces them and is one-to-one on the other chords. Let us say that a chord diagram without pairs of adjacent chords is \emph{reduced}.
\begin{figure}[ht]
$$\includegraphics[height=0.7in]{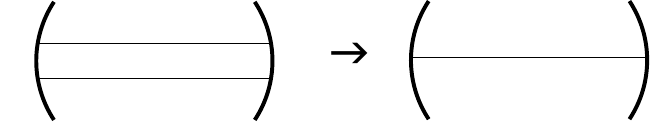}$$
\caption{Reduction of a chord diagram}\label{chordreduction}
\end{figure}

\begin{lemma}\label{reduced}
Let $C'$ be a reduced diagram obtained from a chord diagram $C$ by a sequence of successive reductions. Then
$C'$ does not depend on the order in which the pairs of adjacent chords are chosen.
\end{lemma}

\begin{proof}
If a chord diagram $C$ has two pairs of adjacent chords, these pairs are either disjoint or have one chord in common. In both cases, the result of reducing $C$ with respect to both pairs does not depend on the order. It is obvious in the case when the pairs are disjoint; if they are not, inspection shows that the order of reductions does not matter in this case either. 
\end{proof}

 \emph{Quiver diagrams}\footnote{The term ``quiver'' here has nothing to do with directed graphs and refers to a pack of arrows.} are a generalization of arrow diagrams. 
 A quiver diagram is a chord diagram in which the endpoints of the chords are labelled by nonnegative integers in such a way that the sum of the labels for every chord is at least 1. Quiver diagrams in which every chord has one end labelled with 1 and the other end with 0 may be identified with arrow diagrams if we take the direction of each chord from 0 towards 1. Define the \emph{degree} of a quiver diagram to be the sum of its labels. The \emph{underlying chord diagram} of a quiver diagram is obtained by forgetting the labels of the endpoints.
 
\begin{figure}[ht]
$$\includegraphics[height=0.7in]{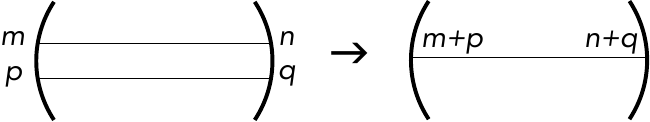}$$
\caption{Reduction of a quiver diagram}\label{quiverreduction}
\end{figure}

One may define the reduction of a quiver diagram with two adjacent chords as the reduction of the underlying chord diagram with the following rule for the labels: the labels of the adjacent endpoints of the adjacent chords are added together while the labels of the other chords remain the same, see Figure~\ref{quiverreduction}. We say that a quiver diagram is reduced if its underlying chord diagram is reduced. As in Lemma~\ref{reduced}, any sequence of reductions of a given quiver diagram always terminates in the same reduced diagram.
Let $\Q$ be the vector space over $\QQ$ spanned by all the quiver diagrams modulo the subspace spanned by the following relations:
\begin{itemize}
\item[(i)] a diagram is equal to its reduction;
\item[(ii)] a diagram with an isolated chord\footnote{that is, a chord whose endpoints are adjacent.} is zero;
\item[(iii)] $\raisebox{-0.17in}{\includegraphics[width=3in]{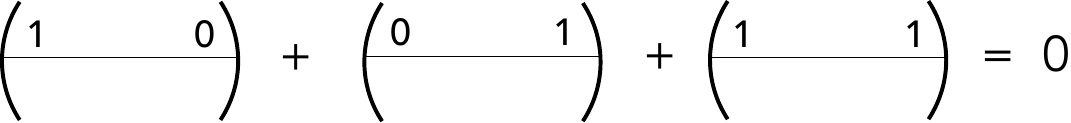}}$ .
\end{itemize} 
Define $\Q_n$ to be the quotient of $\Q$ by the subspace spanned by the diagrams of degrees greater than $n$. For a reduced chord diagram $C$ without isolated chords and at most $n$ chords, let $\Q_n^C\subset \Q_n$ be the subspace spanned by the quiver diagrams whose underlying chord diagram is $C$. 
\begin{proposition}
The vector space $\Q_n$ is the direct sum of all the subspaces $\Q_n^C$. 
\end{proposition}
Indeed, in $\Q_n$ any diagram is equal to its reduction and all the quiver diagrams that participate in the relation (iii) have the same underlying chord diagram.

\begin{proposition}\label{basis}
Given a reduced chord diagram $C$ without isolated chords, mark one endpoint for each chord. The set of all quiver diagrams of degree at most $n$ with the underlying chord diagram $C$ and the labels of the marked endpoints all zero form a basis of $\Q_n^C$.
\end{proposition}
One can think of such quiver diagrams as arrow diagrams whose arrows may have multiplicities.
\begin{proof}
Let us choose one particular chord of a diagram $D$ and assume that in the diagram  
$$D=D(p,q)\ = \ \raisebox{-0.25in}{\includegraphics[width=1in]{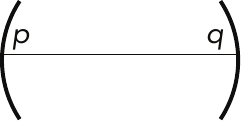}}$$
with $p,q>0$ the marked endpoint of the shown chord is on the right. We have 
\begin{multline*}
 \raisebox{-0.24in}{\includegraphics[width=1in]{diagrampq.pdf}} \quad = \quad  \raisebox{-0.24in}{\includegraphics[width=1in]{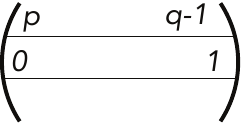}}\quad  =\quad 
 - \quad \raisebox{-0.24in}{\includegraphics[width=1in]{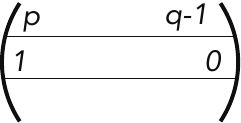}}\quad - \quad \raisebox{-0.24in}{\includegraphics[width=1in]{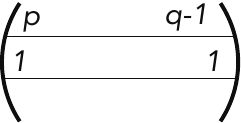}}\\[10pt]
 = \quad - \quad \raisebox{-0.24in}{\includegraphics[width=1in]{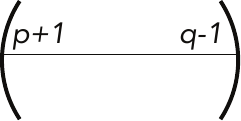}}\quad -\quad  \raisebox{-0.24in}{\includegraphics[width=1in]{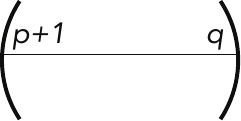}}\quad 
\end{multline*}
or, in other words, $$D(p,q)=-D(p+1,q-1)-D(p+1,q)$$
from where do we get
\begin{equation}\label{redd}D(p,q)=-D(p+1,q-1)+D(p+2,q-1)-D(p+3, q-1)+\ldots.\end{equation}
This sum is finite since any diagram of degree $n+1$ or more is zero in $\Q_n^C$. Now, if $q-1>0$,  the identity (\ref{redd}) can be used to express the diagrams on the right-hand side in terms of the diagrams $D(m, q-2)$ and so on, until we express $D(p,q)$ in terms of diagrams of the form $D(m,0)$.
Then, this algorithm can be applied to all other chords of $D$.

Rewriting in this way the relations in $\Q_n^C$ we obtain a tautology so all the diagrams whose marked endpoints are labelled with 0 are linearly independent.
\end{proof}

Since arrow diagrams can be considered as quiver diagrams and the relations among the diagrams in $\A_n$ also hold in $\Q_n$ we get a map 
$$q_n:\A_n\to \Q_n.$$
The second part of Theorem~\ref{main1} is a consequence of the following statement:
\begin{theorem}\label{main2}
For a non-trivial doodle $D$ with $n$ crossings or fewer, 
$$q_n(I_n(D))\neq 0.$$ 
If $D_1, D_2$ are non-equivalent doodles with $n$ or fewer crossings,   
$$q_{n+1}(I_{n+1}(D_1))\neq q_{n+1}( I_{n+1}(D_2)).$$
\end{theorem}
\begin{figure}[ht]
$$\includegraphics[width=3.5in]{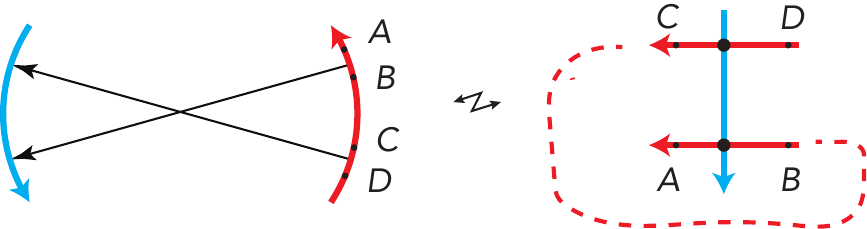}$$
\caption{A pair of intersecting adjacent chords pointing in the same direction and what they might look like in a doodle. The points $B$ and $C$ must be joined in $S^2$ by a segment of the doodle containing no other intersection points, which is impossible.}\label{noparallel2}
\end{figure}
\begin{proof} Assume that the arrow diagram $X_D$ of a minimal doodle $D$ contains a pair of adjacent chords. These chords must point in the same direction since otherwise they could be eliminated by a move $R_2'$; this would mean that the doodle is not minimal. These chords cannot intersect each other, see Figure~\ref{noparallel2}; the only possible situation when a pair of adjacent chords appears in a minimal doodle diagram is when these chords are parallel as shown in Figure~\ref{noparallel3} (a); in fact, in this case there may be a cluster with more than two parallel chords with adjacent endpoints as in Figure~\ref{noparallel3} (b). 
Such a cluster of parallel chords with adjacent endpoints divides the arrow diagram $X_D$ into two parts and no other chord intersects any of the chords of the cluster. The reduction of $X_D$ consists of replacing each maximal cluster of $m$ parallel arrows in $X_D$ with one chord whose initial point is marked by $0$ and the final point by $m$. Denote the resulting reduced quiver diagram by $r(X_D)$, see Figure~\ref{r}. If $X_D$ is reduced to begin with, we set $r(X_D)=X_D$.

\begin{figure}[ht]
$$\includegraphics[width=6.0in]{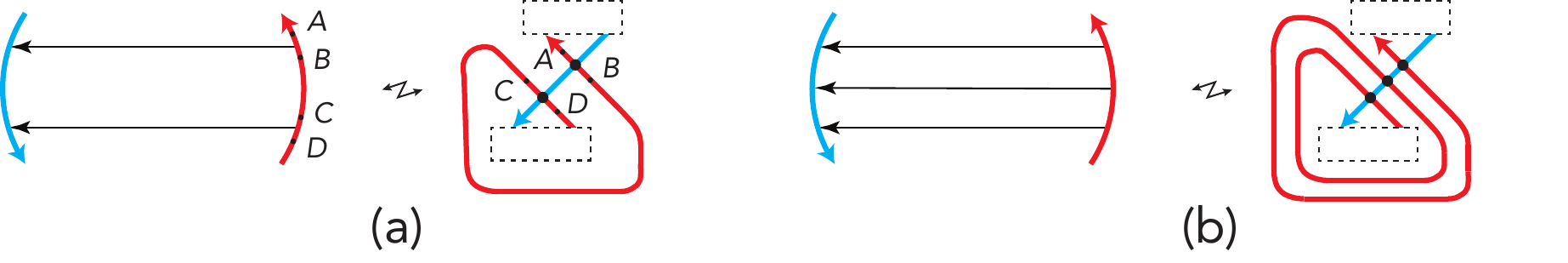}$$
\caption{Doodles whose arrow diagrams have adjacent chords}\label{noparallel3}
\end{figure}

\begin{figure}[ht]
$$\includegraphics[width=3.5in]{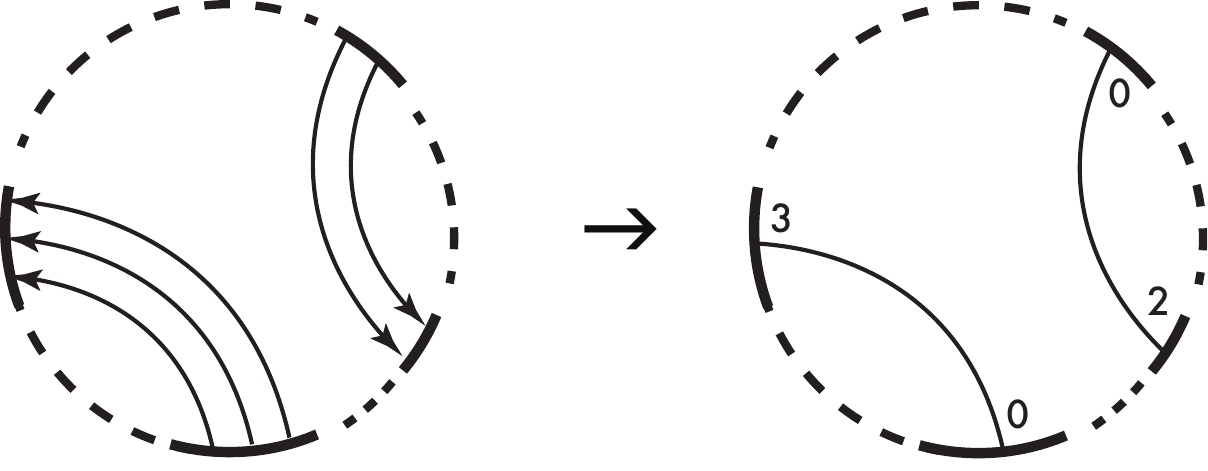}$$
\caption{The diagrams $X_D$ and $r(X_D)$.}\label{r}
\end{figure}

Assume $D$ is minimal and has at most $n$ crossings. Since each chord in $r(X_D)$ has one endpoint marked 0, we can choose a basis in $\Q_n$ as in Proposition~\ref{basis} such that $r(X_D)$ is a basis element. Then, the coefficient of $r(X_D)$ in $q_n(I_n(D))$ is equal to 1 since there is precisely one subdiagram that gives $r(X_D)$; namely, $X_D$ itself. It follows that $q_n(I_n(D))\neq 0$.

Now, let $D_1\neq D_2$ be two minimal doodles with at most $n$ crossings. Without loss of generality, assume that $D_1$ has at least as many crossings as $D_2$. Then, if the underlying chord diagrams of the arrow diagram $X_{D_1}$ and $X_{D_2}$ are different, $r(X_{D_1})\neq r(X_{D_2})$ and we can choose a basis for $\Q_n$ in such a way that both $r(X_{D_1})$ and $r(X_{D_2})$ are basis elements. Then, the coefficient of $r(X_{D_1})$ in $q_n(I_n(D_1))$ is 1 while in $r(X_{D_2})$ it is 0. Indeed, if $r(X_{D_1})$ is a subdiagram of $r(X_{D_2})$, we would have $r(X_{D_1})=r(X_{D_2})$ since the degree of $r(X_{D_1})$ is at least as big as that of $r(X_{D_2})$. This, however, is impossible. 

If the underlying chord diagrams of $X_{D_1}$ and $X_{D_2}$ coincide, $X_{D_2}$ differs from $X_{D_1}$ by the reversal of the directions of some chords. Let $c$ be a chord in $X_{D_2}$ whose direction differs from its direction in $X_{D_1}$. Choose a basis for  $\Q_{n+1}$ so that $X_{D_1}$ is a basis element and let $W$ be the quiver diagram obtained from $X_{D_1}$ by increasing the label of the chord $c$ by 1. Then, the coefficient of $W$ in $q_{n+1}(I_{n+1}(D_1))$ is zero while in $q_{n+1}(I_{n+1}(D_2))$ it is $\pm 1$.

\end{proof}

\section{The diagram invariants are of finite type}\label{Vassiliev}
The theory of finite-type invariants for doodles was developed by Vassiliev in \cite{Vas1}.
A \emph{singular doodle} consists of the following data:
\begin{itemize}
\item an immersed curve $S^1\to S^2$ with a finite number of self-intersection points, such that
 the tangent lines to the branches of the curve at each self-intersection point are all distinct;
\item a linear order in the branches of the curve at each self-intersection point.
\end{itemize}
A \emph{singular point} of a singular doodle is a self-intersection point where three or more branches meet. The \emph{complexity} of a singular point is the number of branches that meet there minus 1. The %\emph{order} 
complexity of a singular doodle is the sum of the multiplicities of its singular points.

Choose a singular point $x$ of a singular doodle $D$ and let $v$ be the tangent vector to the last branch of $D$ at $x$.
Write $v_\perp$ for the vector obtained from $v$ by rotating it by $\pi/2$ in the positive direction. Let $D'$ and $D''$ be the singular doodles obtained by moving the last branch of $D$ slightly off $x$ in the directions $v_\perp$ and $-v_\perp$, respectively. A \emph{one-step resolution} of $D$ at $x$ is the formal difference $D'-D''$. This procedure can be iterated. The order of $D'$ and $D''$ is one less than that of $D$; therefore, by taking consecutive one-step resolutions we arrive at a linear combination of usual doodles that we call the \emph{complete resolution} of $D$. 

An invariant of doodles can be extended to linear combinations of doodles by linearity. An invariant is said to be \emph{of order $m$} if it vanishes on the complete resolution of any singular doodle with complexity greater than $m$.

\begin{proposition}\label{propft}
The $n$th diagram invariant $I_n$ vanishes on any singular doodle of complexity greater than $n+s$ where $s$ is the number of singular points of the singular doodle.
\end{proposition}
Since the complexity of each singular point is at least 2, the number $s$ of singular points is at most half the complexity of the singular doodle. If a singular doodle has complexity $c$ with $c>2n$, then $c>n+s$ and Proposition~\ref{propft} implies that that $I_n$ is zero on it,  
which establishes Theorem~\ref{main1}.  

\begin{proof}[Sketch of the proof of Proposition~\ref{propft}] 
If $\mathcal{R}$ is a region in $S^2$ with piecewise smooth boundary, a \emph{flat tangle}  is an immersion of a 1-manifold into $\mathcal{R}$ which is transversal to the boundary and has a finite number of transversal double points, all in the interior of $\mathcal{R}$. 

Assume $\mathcal{R}_1$ and $\mathcal{R}_2$ intersect along a component of their boundary, and $T$ is a flat tangle in $\mathcal{R}_1 \cup \mathcal{R}_2$ which is transversal to $\partial \mathcal{R}_1$. Then $T$ defines two flat tangles,  $T_1$ in $\mathcal{R}_1$ and $T_2$ in $\mathcal{R}_2$ and we say that $T$ is a product of $T_1$ and $T_2$.

A flat tangle $T$ gives rise to an arrow diagram whose skeleton is the 1-manifold paramet-rizing $T$. One may also talk about the sum of all the subdiagrams of the arrow diagram of $T$, we denote this sum by $I(T)$. In the same way as the singular doodle, we may define singular flat tangles. Proposition~\ref{propft} is a consequence of the following fact:
for the singular flat tangle $T_k$ with one singular point where $k$ branches of $T_k$ intersect, $I(T_k)$ involves only diagrams with at least $k-1$ chords. %See Figure~\ref{singpoint}.

Indeed, this is true for $k=3$:
%\begin{multline*}
% \raisebox{-0.24in}{\includegraphics[width=1in]{diagrampq.pdf}} \quad = \quad  \raisebox{-0.24in}{\includegraphics[width=1in]{dd1.pdf}}\quad  =\quad 
% - \quad \raisebox{-0.24in}{\includegraphics[width=1in]{dd2.pdf}}\quad - \quad \raisebox{-0.24in}{\includegraphics[width=1in]{dd3.pdf}}\\[10pt]
% = \quad - \quad \raisebox{-0.24in}{\includegraphics[width=1in]{dd4.pdf}}\quad -\quad  \raisebox{-0.24in}{\includegraphics[width=1in]{dd5.pdf}}\quad 
%\end{multline*}
\begin{multline*}
I\Bigl(\raisebox{-0.21in}{\includegraphics[width=0.5in]{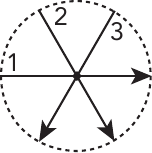}}\Bigr)    
\  = \  I\Bigl( \raisebox{-0.21in}{\includegraphics[width=0.5in]{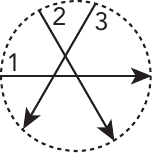}}\Bigr)\     -\ I\Bigl( \raisebox{-0.21in}{\includegraphics[width=0.5in]{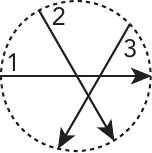}}\Bigr)\\[10pt]
= \raisebox{-0.26in}{\includegraphics[width=0.5in]{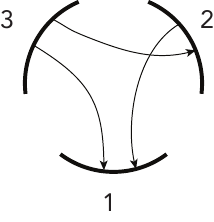}}\,-
\raisebox{-0.26in}{\includegraphics[width=0.5in]{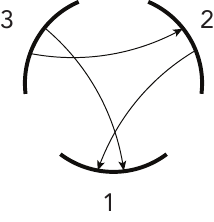}}\,+
\raisebox{-0.26in}{\includegraphics[width=0.5in]{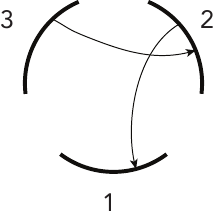}}\,+
\raisebox{-0.26in}{\includegraphics[width=0.5in]{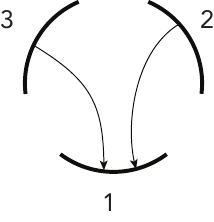}}\,+
\raisebox{-0.26in}{\includegraphics[width=0.5in]{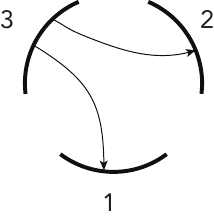}}\,-
\raisebox{-0.26in}{\includegraphics[width=0.5in]{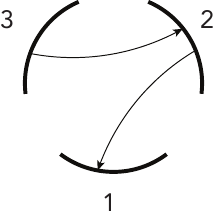}}\,-
\raisebox{-0.26in}{\includegraphics[width=0.5in]{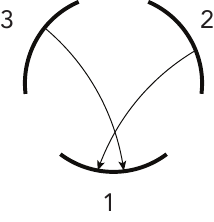}}\,-
\raisebox{-0.26in}{\includegraphics[width=0.5in]{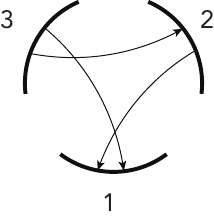}}
\end{multline*}
For $k>3$ use induction. Let $T'-T''$ be a one-step resolution of $T_k$:
$$T'=\raisebox{-0.26in}{\includegraphics[width=0.5in]{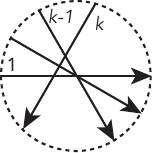}}\, ,\qquad  T''=\raisebox{-0.26in}{\includegraphics[width=0.5in]{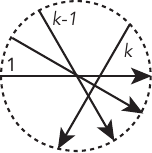}}.$$
Then, $I(T')=I_1(T')+I_2(T')$ where $I_1(T')$ consists of the diagrams that have at least one chord with an endpoint on the $k$th branch and $I_2(T')$ consists of the diagrams that have no such chords. Similarly, we write $I(T'')=I_1(T'')+I_2(T'')$. Now, $I_2(T')=I_2(T'')$. On the other hand, both $I_1(T')$ and $I_2(T'')$ are products of $I(T_{k-1})$ and of a linear combination of tangles with at least one chord. Therefore, by the induction assumption these expressions involve only diagrams with $(k-2)+1=k-1$ chords and the same is true for $$I(T_k)=I_1(T')-I_1(T'').$$
\end{proof}

\section{Some remarks}\label{remx}

\noindent{\sffamily\bfseries 5.1} In the definition of the spaces $\A_n$ and $\Q_n$, the coefficients of the diagrams can be taken in the field $\mathbb{F}_2$ with two elements. The diagram invariants defined in this fashion still distinguish doodles.

\bigskip

\noindent{\sffamily\bfseries 5.2} One may apply the constructions of this note to \emph{long} doodles; these are doodles in $S^2$ that pass through a given point (thought of as the infinity in the Riemann sphere). Long doodles are parametrized by the real line and the preimages of the self-intersection points are naturally ordered. As a consequence, the branches of such a doodle at each self-intersection point are also ordered and each self-intersection point carries a sign which is positive whenever the ordered pair of tangent vectors to the doodle forms a positive basis. Instead of arrow diagrams, long doodles may be encoded by chord diagrams with signed chords; the construction of the subdiagram invariant for such signed chord diagrams is entirely similar to the one described in the present note.

%%% REFERENCES %%%
{\small\bibliography{commat}}
% Please, do not change the above line and do not insert your references
% into this file.  Instead, insert your references into the commat.bib file.
% See commat.bib for further instructions.

\EditInfo{January 19, 2024}{January 27, 2024}{Sergei Chmutov}
\end{document}